\documentclass[12pt]{article}%
\usepackage[applemac]{inputenc}
\usepackage{amsmath,amssymb,fullpage}
\usepackage{showlabels}
\usepackage{amsfonts}
\usepackage{amsmath,amssymb}
\usepackage[applemac]{inputenc}
\usepackage{amsmath,amssymb,fullpage}
\usepackage{color}
\usepackage{color, colortbl}
\usepackage{amsmath}
\usepackage{amssymb}
\usepackage{graphicx}
\usepackage{tikz}
\usepackage{geometry}
\usetikzlibrary{arrows}
\usepackage{amsfonts}
\usepackage{amsmath,amssymb}
\usepackage[applemac]{inputenc}
\usepackage{amsmath,amssymb,fullpage}
\usepackage{color}
\usepackage{amsmath}
\usepackage{amssymb}
\usepackage{graphicx}
\usepackage{tikz}
\newtheorem{theorem}{Theorem}[section]
\newtheorem{lemma}[theorem]{Lemma}
\newcommand{\E}{E}

\newtheorem{definition}[theorem]{Definition}

\newtheorem{example}[theorem]{Example}

\newtheorem{problem}[theorem]{Problem}
\newtheorem{remark}[theorem]{Remark}

\usetikzlibrary{arrows,shapes,automata,petri}

\newcommand{\dproof}{\noindent {Proof.} \quad}
\newcommand{\fproof}{\hfill $\square$ \bigskip}

\numberwithin{equation}{section}

\definecolor{LightCyan}{rgb}{0.88,1,1}

\def\RR{{\mathbb{ R}}}

\def\EE{{\mathbb{ E}}}

\def\1B{\text{1\!\!I}}

\def\l{\langle}
\def\<{\langle}
\def\>{\rangle}
\def\P{\mathbb{P}}
\def\R{\mathbb{R}}
\def\l{\lambda}
\def\S{\mathcal{S}}
\def\E{\mathbb{E}}

\begin{document}

\title{The time-fractional stochastic heat equation driven by time-space white noise}
\author{Rahma Yasmina Moulay Hachemi$^{1}$ \& Bernt \O ksendal$^{2}$
}

\date{31 October 2022}
\maketitle

\footnotetext[1]{%
Laboratory of Stochastic Models, Statistic And Applications University Of Saida- Dr. Moulay Tahar, P.O.Box 138 En-Nasr Saida 20000, Algeria.\newline
Email: yasmin.moulayhachemi@yahoo.com}

\footnotetext[2]{%
Corresponding author\\
Department of Mathematics, University of Oslo, Box 1053 Blindern, 0316 Oslo, Norway.\\ 
Email: oksendal@math.uio.no.}
\paragraph{MSC [2020]:}
\emph{30B50; 34A08; 35D30; 35D35; 35K05; 35R11; 60H15; 60H40}\\

\paragraph{Keywords:}
\emph{Fractional stochastic heat equation; Caputo derivative: Mittag-Leffler function; time-space white noise; tempered distributions: mild solution}

\begin{abstract}
We study the time-fractional stochastic heat equation driven by time-space white noise with space dimension 
  $d\in\mathbb{N}=\{1,2,...\}$ and the fractional time-derivative is the Caputo derivative of order $\alpha \in (0,2)$.   

 We consider the equation  in the sense of distribution, and we find an explicit expression for the $\mathcal{S}'$-valued solution $Y(t,x)$, where $\mathcal{S}'$ is the space of tempered distributions. \\

Following the terminology of Y. Hu \cite{Hu}, we say that the solution is \emph{mild} if $Y(t,x) \in L^2(\P)$ for all $t,x$, where $\P$ is the probability law of the underlying time-space Brownian motion. 
It is well-known that in the classical case with $\alpha = 1$, the solution is mild if and only if the space dimension $d=1$. We prove that if $\alpha \in (1,2)$ the solution is mild if $d=1$ or $d=2$.
 If $\alpha < 1$ we prove that the solution is not mild for any $d$. 
 \end{abstract}

\section{Introduction}
The  \emph{fractional derivative} of a function was first introduced by Niels Henrik Abel in 1823 \cite{A}, in connection with his solution of the tautochrone (isochrone) problem in mechanics.
\vskip 0.2cm
The \emph{Mittag-Leffler function} $E_{\alpha}(z)$ was introduced by  Gösta Magnus Mittag-Leffler in 1903 \cite{ML}. He showed that this function has a connection to the fractional derivative introduced by Abel.
\vskip 0.2cm
The fractional derivative turns out to be useful in many situations, e.g. in the study of waves,  including ocean waves around an oil platform in the North Sea, and ultrasound in bodies. In particular, the fractional heat equation may be used to describe anomalous heat diffusion, and it is related to power law attenuation. This and many other applications of fractional derivatives can be found in the book by S. Holm \cite{H}.

In this paper we study the following fractional stochastic heat equation
\begin{equation}\label{heat}
 	 \frac{\partial^{\alpha}}{\partial t^{\alpha}}Y(t,x)=\l \Delta Y(t,x)+\sigma W(t,x);\; (t,x)\in (0,\infty)\times \mathbb{R}^{d}
 	 	\end{equation}
  	where $d\in\mathbb{N}=\{1,2,...\}$ and $\frac{\partial^{\alpha}}{\partial t^{\alpha}}$ is the Caputo derivative of order $\alpha \in (0,2)$, and $\l>0$ and $\sigma\in \mathbb{R}$ are given constants,
   
  	\begin{equation}
  		\Delta Y =\sum_{j=1}^{d}\frac{\partial^{2}Y}{\partial x_{j}^{2}}(t,x)
  	\end{equation}
  is the Laplacian operator and
  
  \begin{equation}
  	W(t,x)=W(t,x,\omega)=\frac{\partial}{\partial t}\frac{\partial^{d}B(t,x)}{\partial x_{1}...\partial x_{d}}
  \end{equation}
is time-space white noise, $$B(t,x)=B(t,x,\omega); t\geq 0, x \in \R^d, \omega \in \Omega$$
 is time-space Brownian sheet with probability law $\P$.
The boundary conditions are
\begin{align}
    Y(0,x)&=\delta_0(x)\text{ (the point mass at  } 0), \label{1.4}\\
    \lim_{x \rightarrow +/- \text{ }\infty}Y(t,x)&=0.\label{1.5}
\end{align}

In the classical case, when $\alpha=1$, this equation models the normal diffusion of heat in a random  or noisy medium, the noise being represented by the time-space white noise $W(t,x)$. \\
\vskip 0.2cm
 - When $\alpha >1$ the equation models \emph{superdiffusion or enhanced diffusion}, where the particles spread faster than in regular diffusion. This occurs for example in some biological systems.\\
\vskip 0.2cm 
 - When $\alpha <1$ the equation models \emph{subdiffusion}, in which travel times of the particles are longer than in the standard case. Such situation may occur in transport systems.\\

 We consider the equation \eqref{heat} in the sense of distribution, and in  Theorem 1 we find an explicit expression for the $\mathcal{S}'$-valued solution $Y(t,x)$, where $\mathcal{S}'$ is the space of tempered distributions. \\
\vskip 0.2cm
Following the terminology of Y. Hu \cite{Hu}, we say that the solution is \emph{mild} if $Y(t,x) \in L^2(\P)$ for all $t,x$. 
It is well-known that in the classical case with $\alpha = 1$, the solution is mild if and only if the space dimension $d=1$. See e.g.Y. Hu \cite{Hu}. \\
We show that if $\alpha \in (1,2)$ the solution is mild if $d=1$ or $d=2$. \\
Then we show that if $\alpha < 1$ then the solution is not mild for any space dimension $d$.\\

There are many papers dealing with various forms of stochastic fractional differential equations. Some papers which are related to ours are:  \\
 -  In the paper by Kochubel et al \cite{KKd} the fractional heat equation corresponding to random time change in Brownian motion is studied.\\
 - The papers by Bock et al \cite{BGO}, \cite{Bd} 
 are considering stochastic equations driven by grey Brownian motion. \\
  - The paper Röckner et al \cite{LRd} proves existence and uniqueness of general time-fractional linear evolution equations in the Gelfand triple setting.\\
  - The paper  which is closest to our paper is Chen et al \cite{CHN}, where 
a comprehensive discussion is given of a general fractional stochastic heat equations with \emph{multiplicative} noise, and with fractional derivatives in both time and space, is given. In that paper the authors prove existence and uniqueness results as well as regularity results of the solution, and they give sufficient conditions on the coefficients and the space dimension $d$, for the solution to be a random field. \\
Our paper, however,  is dealing with \emph{additive} noise and a more special class of fractional heat equations. As in \cite{CHN} we find explicit solution formulae in the sense of distributions and give conditions under which the solution is a random field in $L^2(\P)$.\\

We refer to Holm \cite{H}, Ibe \cite{I}, Kilbas et al \cite{Kilbas} and Samko et al \cite{samko} for more information about fractional calculus and their applications.

\section{The space of tempered distributions}
For the convenience of the reader we recall some of the basic properties of the Schwartz space $\mathcal{S}$ of rapidly decreasing smooth functions and its dual, the space $\mathcal{S}'$ of tempered distributions.
\subsection{The space of tempered distributions}
Let $n$ be a given natural number. Let $\mathcal{S}=\mathcal{S}(\mathbb{R}^n)$\label{simb-028} be the
space of rapidly decreasing smooth real
functions $f$
on $\mathbb{R}^n$
equipped with the family of seminorms:\label{simb-029} 
\begin{equation*}
\Vert f \Vert_{k,\alpha} := \sup_{y \in \mathbb{R}^n}\big\{ (1+|y|^k) \vert
\partial^\alpha f(y)\vert \big\}< \infty,
\end{equation*}
where $k = 0,1,...$, $\alpha=(\alpha_1,...,\alpha_n)$ is a multi-index with $%
\alpha_j= 0,1,...$ $(j=1,...,n)$ and\label{simb-030} 
\begin{equation*}
\partial^\alpha f := \frac{\partial^{|\alpha|}}{\partial
y_1^{\alpha_1}\cdots \partial y_n^{\alpha_n}}f
\end{equation*}
for $|\alpha|=\alpha_1+ ... +\alpha_n$.

Then
$\mathcal{S}=\mathcal{S}(\mathbb{R}^n)$ is a
Fr\'echet space.

Let $\mathcal{S}^{\prime }=\mathcal{S}^{\prime }(\mathbb{R}^{n})$\label%
{simb-031} be its dual, called the space of \emph{tempered distributions}. 
\index{tempered distributions} Let $\mathcal{B}$ denote the family of all
Borel subsets of $\mathcal{S}^{\prime }(\mathbb{R}^{n})$ equipped with the
weak* topology. If $\Phi \in \mathcal{S}^{\prime }$ and $f \in \mathcal{%
S}$ we let \label{simb-033} 
\begin{equation}
\Phi (f) \text{ or } \langle \Phi ,f \rangle  \label{3.1}
\end{equation}%
denote the action of $\Phi$ on $f$.

\begin{example}
\begin{itemize}
\item
{(Evaluations)}
For $y \in \R$ define the function $\delta_y$ on $\S(\R)$ by $\delta_y(\phi)=\phi(y)$. Then $\delta_y$ is a tempered distribution.\\
\item
{(Derivatives)} Consider the function $D$, defined for $\phi \in \S(\mathbb{R})$ by $D[\phi]=\phi^{\prime}(y)$. Then  $D$ is a tempered distribution. \\
\item
{(Distributional derivative)}\\
 Let $T$ be a tempered distribution, i.e. $T \in \S^{'}(\mathbb{R}) $. We define the distributional derivative $T^{'}$ of $T$ by
 $$ T^{'}[\phi]=-T[\phi^{'}]; \quad \phi \in \S.$$
 Then $T^{'}$ is again a tempered distribution.
 \end{itemize}
 \end{example}
 
 In the following we will apply this to the case when $n=1+d$ and $y=(t,x) \in \R \times \R^d$.
 \section{The Mittag-Leffler functions}
 \begin{definition}(The two-parameter Mittag-Leffler function)
 The Mittag-Leffler function of two parameters $\alpha,\; \beta$ is denoted by $E_{\alpha,\beta}(z)$ and defined by:
\begin{equation}
    E_{\alpha,\beta}(z)=\sum_{k=0}^{\infty}\frac{z^{k}}{\Gamma(\alpha k+\beta)}
\end{equation}
where $z,\; \alpha,\; \beta\in \mathbb{C},\; Re(\alpha)>0\; and\; Re(\beta)>0,$ and $\Gamma$ is the Gamma function.
\end{definition}
 
 \begin{definition}(The one-parameter Mittag-Leffler function)
The Mittag-Leffler function of one parameter $\alpha$ is denoted by $E_{\alpha}(z)$ and defined as;
\begin{equation}
    E_{\alpha}(z)=\sum_{k=0}^{\infty}\frac{z^{k}}{\Gamma(\alpha k+1)}
\end{equation}
where $z,\; \alpha\in \mathbb{C},\; Re(\alpha)>0.$
\end{definition}

\begin{remark} Note that $E_{\alpha}(z)= E_{\alpha,1}(z)$ and that 
\begin{align}
E_{1}(z)=\sum_{k=0}^{\infty}\frac{z^{k}}{\Gamma( k+1)}=\sum_{k=0}^{\infty}\frac{z^{k}}{k!} =e^{z}.
\end{align}
\end{remark}

\section{The (Abel-)Caputo fractional derivative}
  In this section we present the definitions and some properties of the Caputo derivatives.
  \begin{definition}
  The (Abel-)Caputo fractional derivative of order $\alpha > 0$ of a function $f$ such that $f(x)=0$ when $x<0$ is denoted by  $D^{\alpha} f (x)$ or $\frac{d^{\alpha}}{dx^{\alpha}} f(x)$ and defined  by 
  \begin{align}\label{caputo1}
  D^{\alpha}f(x):& =
  \begin{cases}
  \frac{1}{\Gamma(n-\alpha)}\int_0^x \frac{f^{(n)}(u)du}{(x-u)^{\alpha +1 -n}}; \quad n-1 < \alpha < n\\
  \frac{d^n}{dx^n}f(x); \quad \alpha =n.
  \end{cases}
  \end{align}
  Here $n$ is an smallest integer greater than or equal to $\alpha$.\\
  
  \noindent If $f$ is not smooth these derivatives are interpreted in the sense of distributions.
  \end{definition}
  
  \begin{example}
  If $f(x)=x$ and $\alpha \in (0,1)$ then
 
  \begin{align}
 D^{\alpha}f(x)=\frac{ x^{1-\alpha} }{(1-\alpha)\Gamma(1-\alpha)}.
 \end{align} 
 In particular, choosing $\alpha=\tfrac{1}{2}$ we get
 \begin{align}
 D^{\tfrac{1}{2}}f(x)=\frac{2 \sqrt{x}}{\sqrt{\pi}}.
 \end{align} 
  \end{example}

\subsection{Laplace transform of Caputo derivatives} 
Recall that the Laplace transform $L$ is defined by 
 
\begin{equation}
	Lf(s)=\int_{0}^{\infty}e^{-st}f(t)dt=:\widetilde{f}(s)
	\end{equation}
	for all $f$ such that the integral converges.

 Some of the properties of the Laplace transform that we will need are:

  \begin{align}
    &L[ \frac{\partial ^{\alpha}}{\partial t^{\alpha}}f(t)](s)=s^{\alpha}(L f)(s)-s^{\alpha-1}f(0) \label{L1}\\       
     &L[E_{\alpha}(bx^{\alpha})](s)  = \frac{s^{\alpha -1}}{s^{\alpha}-b}\label{L2}\\
     &L[x^{\alpha-1}E_{\alpha,\alpha}(-b x^{\alpha})](s)=\frac{1}{s^{\alpha}+b}\label{L3}
     \end{align}
     Recall that the convolution $f\ast g$ of two functions $f,g: [0,\infty) \mapsto \mathbb{R}$  is defined by
\begin{align}
(f \ast g)(t)=\int_0^t f(t-r)g(r) dr; \quad t \geq 0.
\end{align}

The convolution rule for Laplace transform states that $$L\left( \int_{0}^{t}f(t-r)g(r)dr\right) (s)=Lf(s)Lg(s),$$ 
or 
\begin{equation}\label{12}
	\int_{0}^{t}f(t-w)g(w)dw=L^{-1}\left( Lf(s)Lg(s)\right) (t).
\end{equation}

 \section{Time-space white noise}
Let $n$ be a fixed natural number. Later we will set $n= 1 + d$.
Define $\Omega={\cal S}'(\RR^n)$, equipped
with the weak-star topology. This space will be the base of our basic
probability space, which we explain in the following:
\vskip 0.3cm
As events we will use the family $\mathcal{F}={\cal
B}({\cal S}'(\RR^n))$ of Borel subsets of ${\cal S}'(\RR^d)$, and our
probability measure $\P$ is defined by the following result:

\begin{theorem}{\bf (The Bochner--Minlos theorem)}\\
There exists a unique probability measure $\P$ on ${\cal B}({\cal S}'(\RR^n))$
with the following property:
$$\E[e^{i\langle\cdot,\phi\rangle}]:=\int\limits_{\cal S'}e^{i\langle\omega,
\phi\rangle}d\mu(\omega)=e^{-\tfrac{1}{2} \Vert\phi\Vert^2};\quad i=\sqrt{-1}$$
for all $\phi\in{\cal S}(\RR^n)$, where
$\Vert\phi\Vert^2=\Vert\phi\Vert^2_{L^2(\RR^n)},\quad\langle\omega,\phi\rangle=
\omega(\phi)$ is the action of $\omega\in{\cal S}'(\RR^n)$ on
$\phi\in{\cal S}(\RR^n)$ and $\E=\E_{\P}$ denotes the expectation
with respect to  $\P$.
\end{theorem} 
We will call the triplet $({\cal S}'(\RR^n),{\cal
B}({\cal S}'(\RR^n)),\P)$ the {\it  white noise probability
space\/}, and $\P$ is called the {\it white noise probability measure}.

The measure $\P$ is also often called the (normalised) {\it Gaussian
measure\/} on ${\cal S}'(\RR^n)$. It is not difficult to prove that if $\phi\in L^2(\RR^n)$ and
we choose $\phi_k\in{\cal S}(\RR^n)$ such that $\phi_k\to\phi$ in $L^2(\RR^n)$,
then
$$\langle\omega,\phi\rangle:=\lim\limits_{k\to\infty}\langle\omega,\phi_k\rangle
\quad\text{exists in}\quad L^2(\P)$$
and is independent of the choice of $\{\phi_k\}$. In
particular, if we define
$$\widetilde{B}(x):=\widetilde{B}(x_1,\cdots,x_n,\omega)=\langle\omega,\chi_
{[0,x_1]\times\cdots\times[0,x_n]}\rangle; \quad  
x=(x_1,\cdots,x_n)\in\RR^n,$$
where $[0,x_i]$ is interpreted as $[x_i,0]$ if $x_i<0$,
then $\widetilde{B}(x,\omega)$ has an $x$-continuous version $B(x,\omega)$, which
becomes an \emph{$n$-parameter Brownian motion}, in the following sense:

By an \emph{$n$-parameter Brownian motion} we mean a family
$\{B(x,\cdot)\}_{x\in\RR^n}$ of random variables on a probability space
$(\Omega,{\cal F},\P)$ such that
\begin{itemize} 
\item
$B(0,\cdot)=0\quad\text{almost surely with respect to } \P,$
\item
$\{B(x,\omega)\}$ is a continuous and Gaussian stochastic process 
\item
For all $x=(x_1,\cdots,x_n)$,
$y=(y_1,\cdots,y_n)\in\RR_+^n$,
$B(x,\cdot),\,B(y,\cdot)$ have the covariance
$\prod_{i=1}^n x_i\wedge y_i$. For general $x,y\in\RR^n$ the covariance is
$\prod_{i=1}^n\int_\RR\theta_{x_i}(s)\theta_{y_i}(s)ds$, where
$\theta_x(t_1,\dots,t_n)=\theta_{x_1}(t_1)\cdots\theta_{x_n}(t_n)$, with
\begin{equation*} 
\theta_{x_j}(s)=
\begin{cases}
1\quad \text{ if } 0<s\leq x_j\\
-1\quad \text{ if } x_j<s\leq 0\\
0 \quad \text{ otherwise}
\end{cases}
\end{equation*} 
\end{itemize} 

It can be proved that the process $\widetilde{B}(x,\omega)$ defined above has a modification $B(x,\omega)$ which satisfies all these properties.
This process $B(x,\omega)$ then becomes an \emph{$n$-parameter Brownian motion}. 
\vskip0.3cm
We remark that for $n=1$ we get the classical (1-parameter) Brownian motion
$B(t)$ if we restrict ourselves to $t\geq 0$. For $n \geq 2$ we get what is often
called the \emph{Brownian sheet}.
\vskip0.3cm
With this definition of Brownian motion it is natural to define the
$n$-parameter Wiener--It\^{o} integral of $\phi\in L^2(\RR^n)$ by
$$\int\limits_{\RR^n}\phi(x)dB(x,\omega):=\langle\omega,\phi\rangle;\quad
\omega\in{\cal S}'(\RR^d).$$
We see that by using the Bochner--Minlos theorem we have obtained an easy construction of
$n$-parameter Brownian motion that works for any parameter dimension $n$. Moreover, we get a representation of the space $\Omega$ as the Fr\' echet space $\mathcal{S}'(\R^d)$. This is an advantage in many situations, for example in the construction of the Hida-Malliavin derivative, which can be regarded as a stochastic gradient on $\Omega$.
See e.g. \cite{DOP} and the references therein.

In the following we put $n=1+d$ and let 
 $$B(t,x)=B(t,x,\omega); t \geq 0, x \in \R^d, \omega \in \Omega$$ 
 denote the (1-dimensional) time-space Brownian motion (also called the Brownian sheet) with probaility law $\P$. Since this process is $(t,x)$-continuous a.s., we can for a.a. $\omega \in \Omega$ define its derivatives with respect to $t$ and $x$ in the sense of distributions.
 Thus we define the time-space white noise $W(t,x)=W(t,x,\omega)$ by
 
\begin{equation}
  	W(t,x)=\frac{\partial}{\partial t}\frac{\partial^{d}B(t,x)}{\partial x_{1}...\partial x_{d}}. \label{WN}
  \end{equation}
In particular, for $d=1$ and $x_1=t$ and get the familiar identity
$$W(t)={{d}\over{dt}}B(t)\hbox{ in } \S^{'}.$$ 

The definition \eqref{WN} can also be interpreted as an element of the Hida space $(\S)^*$ of \emph{stochastic distributions}, and in that setting it has been proved (see Lindstr\o m, \O . , Ub\o e \cite{LOU} and Benth \cite{B}) that the Ito-Skorohod integral with respect to $B(dt,dx)$ can be expressed as 
\begin{align}
\int_0^T \int_{\mathbb{R}^d} f(t,x,\omega) B(dt, dx)=\int_0^T \int_{\mathbb{R}^d} f(t,x,\omega) \diamond W(t, x) dt dx,
\end{align}
where $\diamond$ denotes the Wick product.\\

In particular, if $f(t, x,\omega)=f(t, x)$ is deterministic, this gives
\begin{align}
\int_0^T \int_{\mathbb{R}^d} f(t, x) B(dt, dx)=\int_0^T \int_{\mathbb{R}^d} f(t,x)  W(t,x) dt dx.
\end{align}
This is the interpretation we are using in this paper.

\section{The solution of the fractional stochastic heat equation}
We now state and prove the first main result of this paper:
\begin{theorem} \label{th1}
The unique solution $Y(t,x) \in \mathcal{S}'$ of the fractional stochastic heat equation \eqref{heat} - \eqref{1.5} is given by
\begin{align}
Y(t,x)&=I_1 + I_2,
\end{align}
where
\small{
\begin{align}
I_1=(2\pi)^{-d} \int_{\mathbb{R}^d} e^{ixy} E_{\alpha}(- \l t^{\alpha} |y|^2) dy
=(2\pi)^{-d} \int_{\mathbb{R}^d} e^{ixy}\sum_{k=0}^{\infty} \frac{(- \l t^{\alpha} |y|^2)^k}{\Gamma(\alpha k +1)}dy, 
\end{align}
}
and

\small{
\begin{align}
    &I_2= \sigma (2\pi)^{-d} \int_{0}^{t}(t-r)^{\alpha -1}\int_{\mathbb{R}^{d}}\left(\int_{\mathbb{R}^{d}}e^{i(x-z)y} E_{\alpha,\alpha}(-\l (t-r)^{\alpha}|y|^2) dy\right) B(dr,dz)\nonumber\\
    &=\sigma (2\pi)^{-d} \int_{0}^{t}(t-r)^{\alpha -1}\int_{\mathbb{R}^{d}}\left(\int_{\mathbb{R}^{d}}e^{i(x-z)y}\sum_{k=0}^{\infty}\frac{(-\l (t-r)^{\alpha}|y|^2)^{k}}{\Gamma(\alpha k+\alpha))}dy\right) B(dr,dz)  
\end{align}
}
where $|y|^2=y^2=\sum_{j=1}^d y_j^2.$ 
\end{theorem}

\dproof
a) First assume that $Y(t,x)$ is a solution of \eqref{heat}.  We apply the Laplace transform $L$ 
to both sides of  \eqref{heat} and obtain (see \eqref{L1}):
   
\begin{equation}
	s^{\alpha}\widetilde{Y}(s,x)-s^{\alpha-1}Y(0,x)=\lambda\widetilde{\Delta Y}(s,x)+\sigma \widetilde{W}(s,x).
	\end{equation}
Applying the Fourier transform $F$, defined by
  
\begin{equation}
	Fg(y)=\int_{\mathbb{R}}e^{-ixy}g(x)dx=:\widehat{g}(y);\; g\in L^{1}(\mathbb{R}^d),
	\end{equation}
  	 we get, since $\widehat{Y}(0,y)=1$, 
  	 \begin{equation}
  	 	s^{\alpha}\widehat{\widetilde{Y}}(s,y)-s^{\alpha-1}=\l\sum_{j=1}^{d}y_{j}^{2}\widehat{\widetilde{Y}}(s,y)+\sigma\widehat{\widetilde W}(s,y),
  	 \end{equation}
   or, 
   \begin{equation}
   	\left( s^{\alpha}+\l|y|^2\right)\widehat{\widetilde{Y}}(s,y)=  s^{\alpha-1}\widehat{Y}(0^{+},y)+\sigma\widehat{\widetilde{W}}(s,y).
   \end{equation}
   Hence
   \begin{align}
       \widehat{\widetilde{Y}}(s,y)=\frac{s^{\alpha -1}}{s^{\alpha} + \l |y|^2}
       + \frac{\sigma \widehat{\widetilde{W}}(s,y)}{s^{\alpha}+\l |y|^2}.
   \end{align}

   Since the Laplace transform and the Fourier transform commute, this can be written
   \begin{align}
       \widetilde{\widehat{Y}}(s,y)=\frac{s^{\alpha -1}}{s^{\alpha} + \l |y|^2}
       + \frac{\sigma \widetilde{\widehat{W}}(s,y)}{s^{\alpha}+\l |y|^2}.
   \end{align}
 Applying the inverse Laplace operator $L^{-1}$  to this equation we get
 \begin{align}
       \widehat{Y}(t,y)&=L^{-1} \Big(\frac{s^{\alpha -1}}{s^{\alpha} + \l |y|^2}\Big)(t,y)
       + L^{-1}\Big(\frac{\sigma \widetilde{\widehat{W}}(s,y)}{s^{\alpha}+\l |y|^2}\Big)(t,y)\nonumber\\
       &=E_{\alpha,1}(-\l |y|^2 t^{\alpha}) + L^{-1}\Big(\frac{\sigma \widetilde{\widehat{W}}(s,y)}{s^{\alpha}+\l |y|^2}\Big)(t,y),\label{8}
   \end{align}
   where we recall that
   \begin{align} \label{2.12}
   E_{\alpha,\beta}(z)=\sum_{k=0}^{\infty} \frac{z^{k}}{\Gamma(\alpha k + \beta)}    
   \end{align}
   is the Mittag-Leffler function.  
 
   It remains to find 
$L^{-1}\left(\frac{\sigma\widehat{\widetilde{W}}(s,y)}{s^{\alpha}+\l |y|^2}\right)$:\\
Recall that the convolution $f\ast g$ of two functions $f,g: [0,\infty) \mapsto \mathbb{R}$  is defined by
\begin{align}
(f \ast g)(t)=\int_0^t f(t-r)g(r) dr; \quad t \geq 0.
\end{align}
The convolution rule for Laplace transform states that $$L\left( \int_{0}^{t}f(t-r)g(r)dr\right) (s)=Lf(s)Lg(s),$$ 
or 
\begin{equation}\label{12}
	\int_{0}^{t}f(t-w)g(w)dw=L^{-1}\left( Lf(s)Lg(s)\right) (t).
\end{equation}

By \eqref{L3} we have 
\begin{align}
	L^{-1}\left(\frac{1}{s^{\alpha}+\l |y|^2} \right) (t)&=t^{\alpha-1}E_{\alpha,\alpha}(-\l t^{\alpha}|y|^2)\nonumber\\
	&=\sum_{k=0}^{\infty}\frac{t^{\alpha-1}(-\l t^{\alpha}|y|^2)^{k}}{\Gamma(\alpha k+\alpha)}\nonumber\\
	&=\sum_{k=0}^{\infty}\frac{(-\l |y|^2)^{k}t^{\alpha(k+1)-1}}{\Gamma(\alpha(k+1))}\nonumber\\
	&=\sum_{k=0}^{\infty}\frac{(-\l t^{\alpha}|y|^2)^{k}t^{\alpha -1}}{\Gamma(\alpha (k+1))}	\nonumber\\
	&=: \Lambda(t,y).\label{Lambda}
\end{align}
In other words,
\begin{equation}
	\frac{\sigma}{s^{\alpha}+\l |y|^2}=\sigma L \Lambda(s,y),
\end{equation}

Combining with \eqref{12} we get 
\begin{align}
	L^{-1}\left( \frac{\sigma}{s^{\alpha}+\l |y|^2} \widehat{\widetilde{W}} (s,y)\right) (t)&=L^{-1}\left( L\left( \sigma\Lambda(s,y)\right) \widetilde{\widehat{W}}(s,y)\right) (t)\\
	&=\sigma \int_{0}^{t}\Lambda(t-r,y)\widehat{W}(r,y)dr.
\end{align}
Substituting this into \eqref{8} we get 
\begin{equation}
	\widehat{Y}(t,y)=E_{\alpha,1}\left( -\l t^{\alpha}|y|^2\right)+\sigma\int_{0}^{t}\Lambda(t-r,y)\widehat{W}(r,y)dr.
\end{equation}

Taking inverse Fourier transform we end up with
\begin{equation}\label{2.18}
 Y(t,x)=F^{-1}\left( E_{\alpha,1}\left(-\l t^{\alpha}|y|^2 
 \right)\right)(x)+\sigma F^{-1}\left(\int_{0}^{t}\Lambda(t-r,y)\widehat{W}(r,y)dr\right)(x).
\end{equation}
Now we use that  
$$F\left(\int_{\mathbb{R}}f(x-z)g(z)dz\right)(y)=Ff(y)Fg(y),$$
or
\begin{equation}
	\int_{\mathbb{R}}f(x-z)g(z)dz=F^{-1}\Big( Ff(y) Fg(y)\Big)(x).
\end{equation}

This gives 
\begin{align*}
	&F^{-1}\left( \int_{0}^{t}\Lambda(t-r,y)\widehat{W}(r,y)dr\right) (x)\\
	&=\int_{0}^{t}F^{-1}\left( \Lambda(t-r,y)\widehat{W}(r,y)\right) (x)dr\\
	&=\int_{0}^{t}F^{-1}\Big( F\left( F^{-1}\Lambda(t-r,y)\right) (y)FW(r,x)(y)\Big) (x)dr\\
	&=\int_{0}^{t}\int_{\mathbb{R}^{d}}\left( F^{-1}\Lambda(t-r,y)(x-z)\right) W(r,z)dzdr\\
	&=\int_{0}^{t}\int_{\mathbb{R}^{d}}\left( (2\pi)^{-d}\int_{\mathbb{R}^{d}}e^{i(x-z)y}\Lambda(t-r,y)dy\right) W(r,z)dzdr\\
	&=(2\pi)^{-d}\int_{0}^{t}\int_{\mathbb{R}^{d}}\left( \int_{\mathbb{R}^{d}}e^{i(x-z)y}\Lambda(t-r,y)dy\right) B(dr,dz).
\end{align*}

Combining this with \eqref{2.18}, \eqref{2.12} and \eqref{Lambda} we get
\small{
\begin{align*}\label{exp1}
    Y(t,x)&= F^{-1}(\sum_{k=0}^{\infty} \frac{(- \l t^{\alpha} |y|^2)^k}{\Gamma(\alpha k +1)})\\
    &+\sigma (2\pi)^{-d} \int_{0}^{t}\int_{\mathbb{R}^{d}}\left(\int_{\mathbb{R}^{d}}e^{i(x-z)y}\Lambda(t-r,y)dy\right) B(dr,dz)\nonumber\\
    &=(2\pi)^{-d} \int_{\mathbb{R}^d} e^{ixy}\sum_{k=0}^{\infty} \frac{(- \l t^{\alpha} |y|^2)^k}{\Gamma(\alpha k +1)}dy \nonumber\\
    &+ \sigma (2\pi)^{-d} \int_{0}^{t}(t-r)^{\alpha -1}\\
    &\int_{\mathbb{R}^{d}}\left(\int_{\mathbb{R}^{d}}e^{i(x-z)y}\sum_{k=0}^{\infty}\frac{(-\l (t-r)^{\alpha}|y|^2)^{k}}{\Gamma(\alpha (k+1))}dy\right) B(dr,dz). 
\end{align*}
}
This proves uniqueness and also that the unique solution (if it exists) is given by the above formula.
\vskip 0.2cm
b) Next, define $Y(t,x)$ by the above formula. Then we can prove that $Y(t,x)$ satisfies \eqref{heat} by reversing the argument above.
We skip the details.
\fproof

\subsection{The classical case ($\alpha$ = 1)}

It is interesting to compare the above result with the classical case when $\alpha$=1:\\
If $\alpha=1$, we get $Y(t,x)=I_{1}+I_{2}$, where
\begin{equation*}
	I_{1}=(2\pi)^{-d}\int_{\mathbb{R}^{d}}e^{ixy}\sum_{k=0}^{\infty}\frac{\left( -\l t |y|^{2}\right) ^{k}}{k!}dy
\end{equation*}
and
\begin{equation*}
	I_{2}=\sigma(2\pi)^{-d}\int_{0}^{t}\int_{\mathbb{R}^{d}}
	\int_{\mathbb{R}^{d}}e^{i(x-z)y}\sum_{k=0}^{\infty}\frac{ \left( -\l (t-r)|y|^{2}\right) ^{k}}{k!}dyB(dr,dz),
\end{equation*}
where we have used that $\Gamma(k+1)=k!$

By the Taylor expansion of the exponential function, we get
\begin{align*}
 I_{1}&=(2\pi)^{-d}\int_{\mathbb{R}^{d}}e^{ixy}e^{-\l t |y|^{2}} dy\\
 &=(2\pi)^{-d}\left( \frac{\pi}{\l t}\right) ^{\frac{d}{2}}e^{-\frac{|x|^{2}}{4\l t}}\\
 &=(4 \pi \l t) ^{-\frac{d}{2}} e^{-\frac{|x|^{2}}{4\l t}},	
\end{align*}
where we used the general formula
\begin{equation}
\int_{\mathbb{R}^{d}}e^{-\left( a |y|^{2}+2by\right) }dy=\left(\frac{\pi}{a} \right)^{\frac{d}{2}} e^{\frac{b^{2}}{a}};\;a>0;\;b\in\mathbb{C}^d.\label{exp}
\end{equation}

Similarly,
\begin{align*}
	I_{2}&=\sigma(2\pi)^{-d}\int_{0}^{t}\int_{\mathbb{R}^{d}} \int_{\mathbb{R}^{d}}e^{i(x-z)y}\sum_{k=0}^{\infty}\frac{\left(-\l (t-r)|y|^{2} \right)^{k} }{k!}dyB(dr,dz)\\
	&=\sigma(2\pi)^{-d}\int_{0}^{t}\int_{\mathbb{R}^{d}}\left( \frac{\pi}{\l (t-r)}\right) ^{\frac{d}{2}}e^{-\frac{|x-z|^{2}}{4\l (t-r) }}B(dr,dz)\\
	&=\sigma (4 \pi \l)^{-\tfrac{d}{2}}\int_0^t \int_{\mathbb{R}^d} (t-r)^{-\tfrac{d}{2}} e^{-\frac{|x-z|^2}{4 \l (t-r)}} B(dr,dz).
\end{align*}
Summarising the above, we get, for $\alpha=1$, 
\begin{align}
	Y(t,x)&=(4\pi \l t)^{-\frac{d}{2}}e^{-\frac{|x|^{2}}{4\l t}}\nonumber\\
	&+\sigma(4\pi\l)^{-\frac{d}{2}}\int_{0}^{t}\int_{\mathbb{R}^{d}}(t-r)^{-\frac{d}{2}}e^{-\frac{|x-z|^{2}}{ 4\l (t-r)}}B(dr,dz)\;\; \label{alpha1}
\end{align}
This is in agreement with a well-known classical result. See e.g. Section 4.1 in Y.Hu \cite{Hu}. 

\section{When is $Y(t,x)$ a mild solution?}

It was pointed out already in 1984 by John Walsh \cite{W} that (classical) SPDEs driven by time-space white noise $W(t,x); (t,x) \in [0,\infty) \times \R^d$ may have only distribution valued solutions if $d \geq 2$. 
Indeed, the solution $Y(t,x)$ that we found in the previous section is in general distribution valued. But in some cases the solution can be represented as an element of $L^2(\P)$. Following Y. Hu \cite{Hu} we make the following definition:

\begin{definition}
The solution $Y(t,x)$ is called \emph{mild} if $Y(t,x) \in L^2(\P)$ \\for all $t>0, x \in \R^d$.
\end{definition}

The second main issue of this paper is the following:
\begin{problem}
   For what values of $\alpha \in (0,2)$ and what dimensions $d=1,2, ...$ is $Y(t,x)$ mild?
      \end{problem}

Before we discuss this problem,  we prove some auxiliary results:
\begin{lemma}(Abel's test) \label{Abel}\\
Suppose $\sum_{n=1}^{\infty}b_{n}$ is convergent  and put $M= \sup\limits_{n} |b_{n}|$. Let $\left\lbrace \rho_{n}\right\rbrace $  be a bounded monotone sequence, and put $R=\sup\limits_{n}|\rho_{n}|$. Then $\sum_{n=1}^{\infty} b_{n}\rho_{n}$ is convergent, and $|\sum_{n=1}^{\infty} b_{n}\rho_{n}|\leq MR+R|\sum_{n=1}^{\infty} b_{n}|$.
   \end{lemma}
   
\dproof
 By summation by parts we have with\\
 $B_{N}=\sum_{k=1}^{N}b_{k};\;N=1,2,...,$
 \begin{align}
  \sum_{k=1}^{N}b_{k}\rho_{k}&=\sum_{k=0}^{N}\rho_{k}(B_{k}-B_{k-1})\\
  &=\sum_{k=1}^{N-1}B_{k}(\rho_{k}-\rho_{k+1})+\rho_{N}B_{N}.
 \end{align}
 Note that 
 \begin{align}
     |\sum_{k=0}^{N-1}B_{k}(\rho_{k}-\rho_{k+1})|&\leq M|\sum_{k=0}^{N-1}\rho_{k}-\rho_{k+1}|=M(\rho_{1}-\rho_{n})\\
     &\leq M R.
 \end{align}
 
 Hence
 \begin{align}
     |\sum_{k=1}^{N}b_{k}\rho_{k}|\leq MR+R|B_{N}|.
 \end{align}
 \fproof

\begin{lemma}\label{4.4}
Suppose $\alpha > 1$. Define 
\begin{align}
    \rho_k=\frac{\Gamma(k+1)}{\Gamma(\alpha k +1)}; k=1,2, ...
\end{align}
Then $\{\rho_k\}_{k}$ is a decreasing sequence.
\end{lemma}
\dproof\\
Consider
\begin{align}
    \frac{\rho_{k+1}}{\rho_k}&= \frac{\Gamma(k+2) \Gamma(\alpha k+1)}{\Gamma(\alpha (k+1)+1) \Gamma(k+1)}=\frac{(k+1) \Gamma(k+1) \Gamma (\alpha k +1)}{\alpha (k+1) \Gamma(\alpha(k+1)) \Gamma(k+1)}\nonumber\\
    &=\frac{\Gamma(\alpha k +1)}{\Gamma(\alpha k +\alpha)} < 1,
\end{align}
since $\alpha > 1.$
\fproof

\begin{lemma} \label{4.5}
Suppose $\alpha > 1$. Define
\begin{align}
    r_k= \frac{\Gamma(k+1)}{\Gamma(\alpha k+\alpha))}; \quad k=1, 2, ...
\end{align}
Then $\{r_k\}_{k}$ is a decreasing sequence.\\
\end{lemma}
\dproof\\
Consider
\begin{align*}
    \frac{r_{k+1}}{r_k}&= \frac{\Gamma(k+2) \Gamma(\alpha (k+1))}{\Gamma(\alpha(k+2)) \Gamma(k+1)}=\frac{(k+1)\Gamma(k+1)\Gamma(\alpha(k+1))}{(\alpha k + 2\alpha -1) \Gamma(\alpha k + 2 \alpha -1) \Gamma(k+1)} \\
    &=\frac{k+1}{\alpha k + 2 \alpha -1} \cdot \frac{\Gamma(\alpha k + \alpha)}{\Gamma (\alpha k+2\alpha -1)} < 1.
\end{align*}
\fproof

 We now return to the question about mildness: \\
 A partial answer is given in the following:
 
  \begin{theorem} \label{th2}
 Let $Y(t,x)$ be the solution of the $\alpha$-fractional stochastic heat equation.
Then the following holds:
\begin{itemize}
\item
a) If $\alpha = 1$, then $Y(t,x)$ is mild if and only if $d=1$.
\item
b) If $\alpha > 1$ then $Y(t,x)$ is mild if $d=1$ or $d=2$.
\item
c) If $\alpha < 1$ then $Y(t,x)$ is not mild for any  $d.$
\end{itemize}
 \end{theorem} 

 \dproof
 Recall that $Y(t,x)= I_1 +I_2$, with
\begin{align}
&I_1=(2\pi)^{-d} \int_{\mathbb{R}^d} e^{ixy}\sum_{k=0}^{\infty} \frac{(- \l t^{\alpha} |y|^2)^k}{\Gamma(\alpha k +1)}dy, \label{I1}\\\
    &I_2= \sigma (2\pi)^{-d} \int_{0}^{t}(t-r)^{\alpha -1}\int_{\mathbb{R}^{d}}\left(\int_{\mathbb{R}^{d}}e^{i(x-z)y}\sum_{k=0}^{\infty}\frac{(-\l (t-r)^{\alpha}|y|^2)^{k}}{\Gamma(\alpha (k+1))}dy\right) B(dr,dz) \label{I2}.
\end{align}

\textbf{ a) The case $\alpha=1$:}

 This case is well-known, but for the sake of completeness we prove this by our method:\\
 By \eqref{alpha1} and the Ito isometry we get
   \begin{align}
       \EE[Y^2(t,x)]= J_1 + J_2,
   \end{align}
where
\begin{align}
J_1=I_1^2=(4\pi \l t)^{-d} e^{- \frac{\|x \|^2}{2 \l t}}
\end{align}
and, by using \eqref{exp},
\begin{align}
    J_2&= \sigma^2 (4\pi \l)^{-d} \int_0^t (t-r)^{-d}(2\pi \l (t-r))^{\frac{d}{2}} dr\nonumber\\
    &= \sigma^2 2^{-d} (2\pi \l)^{-\frac{d}{2}} \int_0^t (t-r)^{-\frac{d}{2}} dr,
\end{align}
which is finite if and only if $d=1$.\\

\textbf{b) The case $\alpha > 1$}

By the It$\hat{o}$ isometry we have
$\mathbb{E}\left[ Y^{2}\left( t,x\right) \right] =J_{1}+J_{2},$
where
\begin{align}
J_{1}&=(2\pi)^{-2d}\left( \int_{\mathbb{R}^{d}}e^{ixy}\sum_{k=0}^{\infty}\frac{\left( -\l t^{\alpha}|y| ^{2}\right)^{k} }{\Gamma(\alpha k+1)}dy\right) ^{2}\nonumber\\
&=(2\pi)^{-2d}\left( \int_{\mathbb{R}^{d}}e^{ixy}E_{\alpha} (-\lambda t^{\alpha}|y|^2)dy\right) ^{2}
\end{align}
 and
 \small{
 \begin{align}
&J_{2}=\sigma^{2}(2\pi)^{-2d}\int_{0}^{t}\int_{\mathbb{R}^{d}}(t-r)^{2\alpha -2}\left(\int_{\mathbb{R}^{d}}e^{i(x-z)y}\sum_{k=0}^{\infty}\frac{\left(-\l (t-r)^{\alpha}|y|^{2} \right) ^{k}}{\Gamma(\alpha k+\alpha))}dy \right)^{2}dzdr\nonumber\\
&=\sigma^{2}(2\pi)^{-2d}\int_{0}^{t}\int_{\mathbb{R}^{d}}(t-r)^{2\alpha -2}\left(\int_{\mathbb{R}^{d}}e^{i(x-z)y}E_{\alpha,\alpha}(-\lambda (t-r)^{\alpha} |y|^2)dy \right)^{2}dzdr.
\end{align}
}
 
By Abel's test  and Lemma \ref{4.4} and \eqref{exp} we get
\begin{align*}
	J_{1}&=(2\pi)^{-2d}\Big(\int_{\mathbb{R}^{d}}\left( \sum_{k=0}^{\infty}\frac{e^{ixy}\left( -\l t^{\alpha}|y |^{2}\right) ^{k}}{\Gamma(k+1)}\frac{\Gamma(k+1)}{\Gamma(\alpha k+1)} \right) dy\Big)^2\\
	&\leq C_{1}\Big(\int_{\mathbb{R}^{d}}e^{ixy}\sum_{k=0}^{\infty}\frac{\left( -\l t^{\alpha}|y|^{2}\right) ^{k}}{ \Gamma(k+1)}dy\Big)^2\\
	&=C_{1}\Big(\int_{\mathbb{R}^{d}}e^{ixy}e^{-\l t^{\alpha}|y|^{2}}dy\Big)^2\\
	&=C_{1}\left( \frac{\pi}{\l t^{\alpha}}\right) ^{d} e^{- \frac{2|x|^{2}}{\l t^{\alpha}}} < \infty \text{ for all }t>0,x \in \R^d \text{ and for all } d.
\end{align*}
By the Plancherel theorem, Lemma \ref{4.5}  and \eqref{exp} we get

\small{
\begin{align*}
	J_{2}&=\sigma^{2}(2\pi)^{-2d}\int_{0}^{t}(t-r)^{2\alpha-2}  \int_{\mathbb{R}^{d}}\left( \sum_{k=0}^{\infty}\frac{(-\l (t-r)^{\alpha}|x-z|^{2} )^{k}}{\Gamma(\alpha k+\alpha)}\right) ^{2}dzdr \label{J2}\\
	&=\sigma^{2}(2\pi)^{-2d}\int_{0}^{t}\int_{\mathbb{R}^{d}}(t-r)^{2\alpha-2} \nonumber\\
	& \int_{\mathbb{R}^{d}}\left( \sum_{k=0}^{\infty}\frac{(-\l (t-r)^{\alpha}|x-z|^{2} )^{k}}{\Gamma(k+1)}\frac{\Gamma(k+1)}{\Gamma(\alpha k+\alpha)} \right) ^{2}dzdr\nonumber\\
	&\leq C_{2}\int_{0}^{t}(t-r)^{2\alpha-2}\int_{\mathbb{R}^{d}}\left( \sum_{k=0}^{\infty}\frac{(-\l (t-r)^{\alpha}|x-z|^{2} )^{k}}{\Gamma(k+1)}\right) ^{2}dzdr\nonumber\\
	&=C_{2}\int_{0}^{t}(t-r)^{2\alpha-2}\int_{\mathbb{R}^{d}}\left( e^{-\l (t-r)^{\alpha}|x-z|^{2}}\right)^{2}dzdr\nonumber\\
	&=C_{2}\int_{0}^{t}(t-r)^{2\alpha-2}\int_{\mathbb{R}^{d}}\left( e^{-2\l (t-r)^{\alpha}|x-z|^{2}}\right)dzdr\nonumber\\
	&=C_{2}\int_{0}^{t}(t-r)^{2\alpha-2}\left( \frac{\pi}{2\l  (t-r)^{\alpha}}\right) ^{\frac{d}{2}}dr\nonumber\\
	&=C_{3}\int_{0}^{t}(t-r)^{2\alpha-2}(t-r)^{-\frac{\alpha d}{2}}dr\nonumber\\
	&=C_{3}\int_{0}^{t}(t-r)^{2\alpha-2-\frac{\alpha d}{2}}dr.
\end{align*}
}

This is finite if and only if $2\alpha-2-\frac{\alpha d}{2}>-1$, i.e. $d<4-\frac{2}{\alpha}$\\
If $\alpha=1+\epsilon$, then $4-\frac{2}{\alpha}=2+\frac{2\epsilon}{1+\epsilon}>2$ for all 
$\epsilon>0.$\\
Therefore $J_2 < \infty$ for $d=1$ or $d=2$, as claimed.
\fproof
 
\textbf{ c) The case $\alpha < 1$}
\vskip 0.4cm

\dproof
By \eqref{I2} we see that
\small{
\begin{align}
J_2&= \sigma^{2}(2\pi)^{-2d}\int_{0}^{t}(t-r)^{2\alpha-2}  \int_{\mathbb{R}^{d}}\big(E_{\alpha,\alpha}( -\lambda (t-r)^{\alpha} |x-z|^2)\big)^{2}dzdr\nonumber\\
&= \sigma^{2}(2\pi)^{-2d}\int_{0}^{t}(t-r)^{2\alpha-2}  \int_{\mathbb{R}^{d}}\big(E_{\alpha,\alpha}( -\lambda (t-r)^{\alpha} |y|^2)\big)^{2}dydr\nonumber
\end{align}
}

Choose $\beta$ such that $0 < \alpha \leq \beta \leq 1$.\\
A result of Pollard \cite{P}, as extended by Schneider \cite{S}, states that the map
\begin{equation}\label{A}
x\mapsto h(x):=E_{\alpha,\beta}(-x);\; x\in \mathbb{R}^{d}    
\end{equation}
is completely monotone, i.e,
\begin{equation}\label{B}
(-1)^{n}\frac{d^{n}}{dx^{n}}h(x)\geq 0\;for\;all\;n=0,1,2,...;\;x\in \mathbb{R}^{d}.
\end{equation}
Therefore by Bernstein's theorem there exists a positive, $\sigma$-finite measure $\mu$ on $\mathbb{R}^{+}$ such that
\begin{equation}\label{C}
E_{\alpha,\beta}(-x)=\int_{0}^{\infty}e^{-xs}\mu(ds).
\end{equation}

In fact, it is known that $\mu$ is absolutely continuous with respect to Lebesgue measure and 
\begin{equation} \label{D}
t^{\beta-1}E_{\alpha,\beta}(-t^{\alpha})=\int_{0}^{\infty}e^{-st}K_{\alpha,\beta}(s)ds
\end{equation}
where
\begin{equation}\label{E}
 K_{\alpha,\beta}(s)=\frac{s^{\alpha-\beta}\left[sin((\beta-\alpha)\pi)+s^{\alpha}sin(\beta\pi)\right]}{\pi\left[s^{2\alpha}+2s^{\alpha}cos(\alpha\pi)+1\right]}  
\end{equation}
See Capelas de Oliveira et al \cite{CMV}, Section 2.3.

Putting $t^{\alpha}=x$ this can be written

\begin{equation}\label{F}
E_{\alpha,\beta}(-x)=x^{\frac{1-\beta}{\alpha}}\int_{0}^{\infty}e^{-s x^{\frac{1}{\alpha}}}K_{\alpha,\beta}(s)ds;\;x>0.
\end{equation}
This gives
\begin{equation}\label{F}
E_{\alpha,\beta}(-\rho |y|^2)=\rho^{\frac{1-\beta}{\alpha}}|y|^{\frac{2(1-\beta)}{\alpha}}\int_{0}^{\infty}e^{-s \rho^{\frac{1}{\alpha}}|y|^{\frac{2}{\alpha}}}K_{\alpha,\beta}(s)ds.
\end{equation}
It follows that
\begin{align}
\big(E_{\alpha,\beta}(-\rho |y|^2)\big)^2
 &\sim \big(\rho^{\frac{1-\beta}{\alpha}}|y|^{\frac{2(1-\beta)}{\alpha}}
 \rho^{\frac{-1}{\alpha}}|y|^{\frac{-2}{\alpha}}\big)^2\nonumber\\
 &= \rho^{-\frac{2\beta}{\alpha}} |y|^{-\frac{4\beta}{\alpha}}
\end{align}
Hence, by using polar coordinates we see that
\begin{align}
\int_{\R^d} \big(E_{\alpha,\beta}(-\rho |y|^2)\big)^2dy 
\sim \int_0^{\infty} R^{-\frac{4\beta}{\alpha}} R^{d-1} dR =\infty,
\end{align}
for all $d$.\\
Therefore $J_2 = \infty$ for all $d$.
\fproof

\begin{remark}
\begin{itemize}
\item
See Y. Hu \cite{Hu}, Proposition 4.1 for a generalisation of the above result in the case $\alpha=1$.
\item
In the cases $\alpha > 1, d \geq 3$ we do not know if the solution $Y(t,x)$ is mild or not. This is a topic for future research.
\end{itemize}
\end{remark}
 
\section{Examples}
\subsection{Example 1}
Let us consider the following heat equation where $\alpha<1.$ 
In this case our  equation models \emph{subdiffusion}, in which travel times of the particles are longer than in the standard case. Such situation may occur in transport systems. For $\alpha=\frac{1}{2}$ and $d=2$ we get
\begin{equation} 
	\frac{\partial^{\frac{1}{2}}}{\partial t^{\frac{1}{2}}}Y(t,x)=\lambda \Delta Y(t,x)+\sigma W(t,x);\; (t,x)\in (0,\infty)\times \mathbb{R}^{2}
\end{equation}
The solution is given by:
\begin{align}
	Y(t,x)&=I_1 + I_2,
\end{align}
where
\small{
	\begin{align}
		I_1=(2\pi)^{-2} \int_{\mathbb{R}^2} e^{ixy} E_{\frac{1}{2}}(- \l t^{\frac{1}{2}} |y|^2) dy
		=(2\pi)^{-2} \int_{\mathbb{R}^2} e^{ixy} erfc(-\lambda t^{\frac{1}{2}}|y|^{2})^{\frac{1}{2}}dy, 
	\end{align}
}
(with 
$erfc(z)=\frac{2}{sqrt (\pi)}\int_{0}^{z}
exp(-t^{2})dt$) 
and

\small{
	\begin{equation}
		I_2= \sigma (2\pi)^{-2} \int_{0}^{t}(t-r)^{\frac{1}{2} -1}\int_{\mathbb{R}^{2}}\left(\int_{\mathbb{R}^{2}}e^{i(x-z)y} E_{\frac{1}{2},\frac{1}{2}}(-\l (t-r)^{\frac{1}{2}}|y|^2) dy\right) B(dr,dz)
\end{equation}	 
}

By the Theorem \ref{th2} this solution is not mild.
\subsection{Example 2}

Next, let us consider the heat equation for $\alpha=\frac{3}{2}$. In this case the equation models \emph{superdiffusion or enhanced diffusion}, where the particles spread faster than in regular diffusion. This occurs for example in some biological systems. Now the equation gets the form
\begin{equation} 
	\frac{\partial^{\frac{3}{2}}}{\partial t^{\frac{3}{2}}}Y(t,x)=\lambda \Delta Y(t,x)+\sigma W(t,x);\; (t,x)\in (0,\infty)\times \mathbb{R}^{2}
\end{equation}
By Theorem \ref{th1} the solution is 
\begin{align}
	Y(t,x)&=I_1 + I_2,
\end{align}
where
\small{
	\begin{align}
		I_1=(2\pi)^{-2} \int_{\mathbb{R}^2} e^{ixy} E_{\frac{3}{2}}(- \l t^{\frac{3}{2}} |y|^2) dy
		=(2\pi)^{-2} \int_{\mathbb{R}^2} e^{ixy}\sum_{k=0}^{\infty} \frac{(- \l t^{\frac{3}{2}} |y|^2)^k}{\Gamma(\frac{3}{2} k +1)}dy, 
	\end{align}
}
and

\small{
	\begin{align}
		&I_2= \sigma (2\pi)^{-2} \int_{0}^{t}(t-r)^{\frac{3}{2} -1}\int_{\mathbb{R}^{2}}\left(\int_{\mathbb{R}^{2}}e^{i(x-z)y} E_{\frac{3}{2},\frac{3}{2}}(-\l (t-r)^{\frac{3}{2}}|y|^2) dy\right) B(dr,dz)\nonumber\\
		&=\sigma (2\pi)^{-2} \int_{0}^{t}(t-r)^{\frac{1}{2} }\int_{\mathbb{R}^{2}}\left(\int_{\mathbb{R}^{d}}e^{i(x-z)y}\sum_{k=0}^{\infty}\frac{(-\l (t-r)^{\frac{3}{2}}|y|^2)^{k}}{\Gamma(\frac{3}{2} k+\frac{3}{2}))}dy\right) B(dr,dz)  
	\end{align}
}

By Theorem \ref{th2} this solution is  mild.

\section{Acknowledgments}
We are grateful to Wolfgang Bock for helpful comments.

\section{Declarations}
Statement about no conflict of interest and no financial interests.\\

On behalf of my coauthor and myself I declare the following:\\
•	We have no relevant financial or non-financial interests to disclose.\\
•	We have no competing interests to declare that are relevant to the content of this article.\\
•	We certify that we have no affiliations with or involvement in any organization or entity with any financial interest or non-financial interest in the subject matter or materials discussed in this manuscript.\\
•	We have no financial or proprietary interests in any material discussed in this article.\\


\begin{thebibliography}{99}
   
   \bibitem{A} Abel, N. H. (1823): Oppl\o sning av et par oppgaver ved hjelp av bestemte integraler (in Norwegian). Magazin for naturvidenskaberne, pp. 55-68.
   
   \bibitem{B} Benth, F. E. (1993). Integrals in the Hida distribution space
(S)*. In B. Lindstr\o m, B. \O ksendal, and A.S. Üstünel (editors): "Stochastic
Analysis and Related Topics", Vol. 8, 89-99. Gordon \& Breach.

\bibitem{BGO} Bock, W., Grothaus, M. \& Orge, K.:Stochastic analysis for vector valued generalized grey Brownian motion. ArXiv: 2111.09229v1, 17 Nov. 2021.

\bibitem{Bd} Bock, W. \& da Silva, J. L.:Wick type SDEs driven by grey Brownian motion. AIP Conference Proceedings, 1871(1): 020004, 2017.

\bibitem{CMV} Capelas de Oliveira, E., Mainardi,F.Vaz, J. (2011): Models based on Mittag Leffler functions for anomalous relaxation in dielectrics. arXiv: 1106.1761  v2, 13 Feb. 2014.

\bibitem{CGS}Chen, L. , Guo, Y. \& Song, J. (2022): Moments and asymptotics for a class of  SPDEs with space-time white noise. arXiv 2206.10069v1. 

\bibitem{CHN}Chen, L. , Hu, Y. \& Nualart, D. (2019): Nonlinear stochastic time-fractional slow and fast diffusion equations on $\mathbb{R}^d$. Stochastic Processes and their Applications 129, 5073-5112.
  
   \bibitem{Hu} Hu, Yaozhong (2019): Some recent progress on stochastic heat equations. Acta Mathematica Scientia 39B(3); 874-914.
   
   \bibitem{H} Holm, S. (2019): Waves with Power-Law Attenuation. Springer.
        
   \bibitem{I} Ibe, Oliver C. (2013): Markov Processes for Stochastic Modelling. 2$^{nd}$ edition. Elsevier.
 
 \bibitem{Kilbas} Kilbas, A.A., Srivastava, H.M., Trujillo, J.J. (2006):
	 Theory and Applications of Fractional Differential Equations, Elsevier  Science B.V.
	 
	 \bibitem{KKd} Kochubel, A. N., Kondratiev, Y. \& da Silva, J. L. (2021): On fractional heat equation. Fractional Calculus \& Applied Analysis 24 (1), 73-87.
	    
\bibitem{LOU} Lindstr\o m, T., \O ksendal, B. \& Ub\o e, J. (1992):
Wick multiplication and Ito-Skorohod stochastic differential equations. 
In S. Albeverio et al (editors): "Ideas and Methods in Mathematical Analysis, Stochastics and Applications". Cambridge Univ. Press, pp.183-206.

\bibitem{LRd} Liu, W., Röckner, M. \& da Silva, J. L.: Quasilinear (stochastic) partial differential equations with time-fractional derivatives (2018), SIAM J. Math. Anal. 50(3), 2588-2607.
  
    \bibitem{ML} Mittag-Leffler, M.G. (1903): Sur la nouvelle fonction E(x). C. R. Acad. Sci. Paris 137, 554?558
  
   \bibitem{MS} Meerschaert, M. M., Sikorskii, A. (2019): Stochastic Models for Fractional Calculus, $2^{nd}$ edition. De Greuter.
   
\bibitem{P} Pollard, H. (1948): The completely monotone character of the Mittag-Leffler function $E_{\alpha}(-x)$. Bull. Amer. Math-Soc.54, 1115-1116.
    
	\bibitem{Samko} Samko, S.G., Kilbas, A.A., Marichev, O.I.: Fractional Integrals and Derivatives. Theory and Applications. Gordon and Breach Science Publishers, New York (1993).

\bibitem{S} Schneider, W.B. (1996): Completely monotone generalized Mittag Leffler functions. Expositiones Mathematicae 14, 3-16.


\bibitem{W} Walsh, J. (1984): An Introduction to Stochastic Partial Differential Equations. Springer Lecture Notes.

\end{thebibliography}
\end{document}